\documentstyle[12pt,amsfonts,amssymb,leqno]{article}

\newtheorem{proposition}{Proposition}[section]
\newtheorem{thm}[proposition]{Theorem}

\newtheorem{lemma}[proposition]{Lemma}
\newtheorem{defn}[proposition]{Definition}

\begin{document}

\begin{center}
{\Large \bf $P \not= NP$ for infinite time Turing machines}
\end{center}

\begin{center}
\renewcommand{\thefootnote}{\fnsymbol{footnote}}
\renewcommand{\thefootnote}{arabic{footnote}}
\renewcommand{\thefootnote}{\fnsymbol{footnote}}
{\large Ralf Schindler}${}^{a}$\footnote[1]{\noindent
2000 {\it Mathematics Subject Classification.} 68Q15, 68Q17, 03E15.\\
Keywords: computer science/complexity theory/descriptive set theory.}
\renewcommand{\thefootnote}{arabic{footnote}}
\end{center}
\begin{center} 
{\footnotesize
%${}^{a}${\it Department of Mathematics, University of California, Berkeley, CA
%94720, USA}\\
%${}^a${\it Lehrstuhl f\"ur Mathematische Logik, Humboldt Universit\"at,
%Unter den Linden 6, 
%10099 Berlin, Germany}\\
${}^a${\it Institut f\"ur formale Logik, Universit\"at Wien, 1090 Wien, Austria}} 
\end{center}

\begin{center}
{\tt rds@logic.univie.ac.at}

%\end{center}
%\begin{center}
{\tt http://www.logic.univie.ac.at/${}^\sim$rds/}\\
\end{center}
%

%\bigskip
{\small {\bf Abstract.}
We state a version of the ``$P=$?$NP$'' problem for infinite time Turing
machines. It is observed that $P \not= NP$ for this version.}

\bigskip
%\bigskip
%\section{Introduction.}
Infinite time Turing machines were introduced in \cite{joel}. The purpose of this
note is to point out that $P \not= NP$ for a transfinite version of 
the ``$P=$?$NP$'' problem.

%Infinite time Turing machines produce partial function from the reals to the reals. 
A
``real'' is considered an element of $2^\omega$, i.e., an infinite sequence of $0$'s
and $1$'s. We'll write ${\mathbb R}$ for $2^\omega$.
It will be convenient to think of a Turing machine to come with {\em two}
halting states, the accept state, and the reject state. 

\begin{defn}
Let $A \subset {\mathbb R}$. We say that $A$ is decidable in polynomial time, or 
$A \in P$, if there are a Turing machine $T$ and some $m<\omega$ such that 

\noindent (a) $T$ decides $A$ (i.e, $x \in A$ iff $T$ accepts $x$), and

\noindent (b) $T$ halts on all inputs after $< \omega^m$ many steps. 
\end{defn}

With infinite time Turing machines, all inputs (i.e., reals) 
should be counted as having the same
length, namely $\omega$. So it appears reasonable to have a polynomial time Turing
machine to be one which always halts after $< \omega^m$ many steps, for some fixed
$m<\omega$.

\begin{defn}
Let $A \subset {\mathbb R}$, and let $\alpha \leq \omega_1+1$. We say that $A$ is 
in $P_\alpha$ if there are a Turing machine $T$ and some $\beta < \alpha$ such that 

\noindent (a) $T$ decides $A$ (i.e, $x \in A$ iff $T$ accepts $x$), and

\noindent (b) $T$ halts on all inputs after $< \beta$ many steps. 
\end{defn}

Of course, $P = P_{\omega^\omega}$. Moreover, $P_{\omega_1+1}$ is just the class of
all $A \subset {\mathbb R}$ which are decided by some Turing machine. Also:

\begin{lemma} {\it \bf (\cite{joel}, Theorem 2.6)}
Let $A \subset {\mathbb R}$. Then $A \in P_{\omega^2}$ 
if and only if ``$x \in A$'' is an
arithmetic statement, uniformly in $x$ (i.e., if $A$ is in $\Sigma^0_\omega$). 
\end{lemma}

\begin{lemma}\label{ref}
Let $A \subset {\mathbb R}$. Then 
$A \in P_{\omega_1^{\rm \it CK}}$ if and only if 
$A$ is a hyperarithmetic set, and if
$A \in P_{\omega_1}$ then 
$A$ is a Borel set.
\end{lemma}

{\sc Proof.} The first part is \cite{joel} Theorem 2.7, and the second part is an
immediate consequence of the proof thereof. $\square$ 

\bigskip
We now turn to the class $NP$.

\begin{defn}
Let $A \subset {\mathbb R}$. We say that $A$ is verifiable in polynomial time, or 
$A \in NP$, if there are a Turing machine $T$ and some $m<\omega$ such that 

\noindent (a) $x \in A$ if and only if $(\exists y \ $ $T$ accepts $x \oplus y)$, and

\noindent (b) $T$ halts on all inputs after $< \omega^m$ many steps. 
\end{defn}

\begin{defn}
Let $A \subset {\mathbb R}$, and let $\alpha \leq \omega_1+1$. We say that 
$A$ is in 
$NP_\alpha$, if there are a Turing machine $T$ and some $\beta < \alpha$ such that 

\noindent (a) $x \in A$ if and only if $(\exists y \ $ $T$ accepts $x \oplus y)$, and

\noindent (b) $T$ halts on all inputs after $< \beta$ many steps. 
\end{defn}

Again, $NP = NP_{\omega^\omega}$. $NP_\alpha$ is the class of all projections of sets
in $P_\alpha$. Hence every set in $NP_{\omega_1+1}$ is analytic.

It is now immediate that $P \not= NP$.

\begin{thm}
$NP_{\omega+1} \setminus P_{\omega_1} \not= \emptyset$.
\end{thm}

{\sc Proof.} Let $A$ be an appropriate universal $\Sigma^1_1$ set. There is then
a recursive $R$ such
that $$x \in A \Leftrightarrow \exists y \in {\mathbb R} \ \forall n < \omega \
(x \upharpoonright n, y \upharpoonright n) \in R.$$
For given $x$ and $y$, whether or not $\forall n < \omega \
(x \upharpoonright n, y \upharpoonright n) \in R$ can be checked in time $\omega$.
Therefore, $A \in NP_{\omega+1}$. By \ref{ref}, 
$A$ cannot be in $P_{\omega_1}$, as it is not
Borel. $\square$

\end{document}